\documentclass[11pt]{article}
\usepackage{amsmath}
\usepackage{amsthm}
\usepackage{amssymb}

\newtheorem{theo}{Theorem}[section]

\newtheorem{lemma}[theo]{Lemma}

\newcommand{\HH}{{\cal H}}

\begin{document}
\date{}
\title{
On bipartite coverings of graphs and multigraphs
}
\author{Noga Alon
\thanks
{Princeton University,
Princeton, NJ 08544, USA and
Tel Aviv University, Tel Aviv 69978,
Israel.
Email: {\tt nalon@math.princeton.edu}.
Research supported in part by
NSF grant DMS-2154082.}
}
\maketitle
\begin{abstract}
A bipartite covering of a (multi)graph $G$ is a collection of bipartite
graphs, so that each edge of $G$ belongs to at least
one of them.  The capacity of the covering is the sum 
of the numbers of vertices of these bipartite graphs. In this note
we establish a (modest) strengthening of
old results of Hansel and of Katona and Szemer\'edi,
by showing that the capacity of any bipartite covering of a graph
on $n$ vertices in which the maximum size of an independent set containing
vertex number $i$ is $\alpha_i$, is at least 
$\sum_i \log_2 (n/\alpha_i).$ We also obtain slightly improved bounds
for a recent result of Kim and Lee about the minimum possible capacity
of a biparite covering of complete multigraphs.
\end{abstract}

\section{Introduction}
A bipartite covering $\HH=\{H_1, \ldots ,H_m\}$
of a graph $G$ on the set of vertices $[n]=\{1,2, \ldots ,n\}$
is a collection of bipartite
graphs $H_i$ on $[n]$, so that each edge of $G$ belongs to at least
one of them. Note that each $H_i$ is not necessarily a subgraph of $G$,
the only assumption is that it is a bipartite subgraph of the complete
graph on $[n]$.
The capacity $\mbox{cap}(\HH)$ 
of the cover is the sum  $\sum_i|V(H_i)|$
of the numbers of vertices of these bipartite graphs. A known result of
Hansel \cite{Ha} is that the
capacity of any bipartite covering of the complete graph $K_n$
on $n$ vertices is at least
$n \log_2 n$. This bound is tight when $n$ is  power of $2$.

In \cite{KL} Kim and Lee consider the analogous
problem, where the complete graph $K_n$ is replaced by the complete
multigraph $K_n^{\lambda}$ in which every pair of distinct vertices
is connected by $\lambda$ parallel edges. A bipartite covering here
is a collection of bipartite graphs so that each edge belongs to at
least $\lambda$ of them. They prove that the capacity of each
bipartite covering of $K_n^{\lambda}$ is at least
$$
\max\{2 \lambda (n-1), n [\log n+ \lfloor (\lambda-1)/2 
\rfloor \log (\frac{\log n}{\lambda}) -\lambda -1]\},
$$
where all logarithms here and in the rest of this note are in base
$2$.
They also establish an upper bound: there exists a bipartite
covering of  $K_n^{\lambda}$ of capacity at most
$$
n( \log (n-1)+(1+o(1)) \lambda \log \log n).
$$
This shows that for $\lambda =(\log n)^{1-\Omega(1)}$ the smallest possible
capacity is $n \log n+ \Theta(n \lambda \log \log n)$ but leaves 
an additive gap of $\Omega(\lambda n \log \log n)$ between the
upper and lower bounds. 

The proofs in \cite{KL} proceed by studying a more general problem 
regarding graphons. Our first contribution in this note
is a shorter combinatorial proof
of the results above,
slightly improving the bounds.  Let $\mbox{cap}(n, \lambda)$ denote
the minimum possible capacity of a bipartite covering of
$K_n^{\lambda}$.
\begin{theo}[Lower bound]
\label{t11}
For positive integers $n \geq 2$ and 
$\lambda$, 
$$
\mbox{cap}(n, \lambda) \geq 
\max\{2 \lambda (n-1), n [\log n+ \lfloor (\lambda-1)/2 
\rfloor \log (\frac{2\log n}{\lambda-1})]\}
$$
\end{theo}
\begin{theo}[Upper bound]
\label{t12}
Let $k(n,\lambda)$ denote the minimum length of a binary error
correcting code with distance at least $\lambda$ which has 
at least $n$ codewords. Then 
$\mbox{cap}(n,\lambda) 
\leq n \cdot k(n, \lambda)$.  Therefore
\begin{enumerate}
\item
For any $n \geq 2$
$$\mbox{cap}(n,2) \leq n (\lceil \log n \rceil+1)< n (\log n+2).$$
\item
For any $n$ and $\lambda \leq 0.5 \log n$
$$\mbox{cap}(n,\lambda) \leq n[\log n+ 
(\lambda-1) (\log (\frac{\log n}{\lambda-1})+4)]$$
\item
For any $0<c<1/2$, and for
$\lambda \geq c\frac{\log n}{1-H(c)}$ where
$H(x)=-x \log x -(1-x) \log (1-x)$ is the binary entropy function,
$$\mbox{cap}(n,\lambda) \leq \frac{\lambda}{c} n.$$
\item
For any fixed $\lambda$ there are infinitely many values of $n$ so that
$$\mbox{cap}(n,\lambda) \leq 
n[\log n+ \lfloor (\lambda-1)/2) \rfloor \log \log n+2].
$$
\end{enumerate}
\end{theo}

Katona and Szemer\'edi \cite{KS} proved the following generalization
of the result of Hansel, dealing with the capacity of bipartite coverings 
of general graphs.
\begin{theo}[\cite{KS}]
\label{t13}
Let $G$ be a graph on the set of vertices $[n]$ and let 
$d_1, d_2, \ldots ,d_n$ denote the degrees of its vertices. Then the 
capacity of any bipartite covering of $G$ is at least
$$\sum_{i=1}^n \log(\frac{n}{n-d_i} ).$$
\end{theo}

Our second contribution here is the following strengthening of this result.
\begin{theo}
\label{t14}
Let $G$ be a graph on the set of vertices $[n]$. For each vertex
$i$ let $\alpha_i $ denote the maximum size of an independent set
of $G$ that contains the vertex $i$. Then
the capacity of any bipartite covering of $G$ is at least
$$\sum_{i=1}^n \log(\frac{n}{\alpha_i} ).$$
\end{theo}

Since it is clear that $\alpha_i \leq n-d_i$ for every $i$, this is
indeed a strengthening of the Katona-Szemer\'edi result (Theorem
\ref{t13}). The binomial random graph $G=G(n,0.5)$ is one example for which
Theorem \ref{t14} is strictly stronger than Theorem \ref{t13}. Indeed,
with high probability for $G=G(n,0.5)$, $d_i=(1/2+o(1))n$ for every $i$
and $\alpha_i=(2+o(1)) \log n$ for every $i$. Therefore the lower  bound
of Theorem \ref{t13} for this $G$ is typically 
$(1+o(1)) n$, whereas the lower bound provided by Theorem \ref{t14}
is $n \log n -(1+o(1)) n \log \log n$. This is tight 
since the chromatic number of $G=G(n,0.5)$ is, with high probability,
$\chi_n=(1+o(1)) \frac{n}{2 \log n}$ implying that $G$ admits a 
bipartite covering consisting of 
$\lceil \log \chi_n \rceil$ spanning bipartite graphs, and the
corresponding 
capacity is $n \log n -(1+o(1)) n \log \log n$.

The rest of this note contains the (short) proofs of the results above.

\section{Complete multigraphs}

\subsection{The lower bound}
Let $\HH=\{H_1, \ldots ,H_m\}$ be a bipartite covering of
the complete multigraph 
$K_n^{\lambda}$ on the set of vertices $[n] =\{1,2,\ldots ,n\}$.
We first prove that $\mbox{cap}(\HH) \geq 2 \lambda (n-1)$.
Let $n_i$ denote the number of vertices of $H_i$. Since it is
bipartite the number of its edges is at most $n_i^2/4$. As the
edges of all these graphs cover each of the $n(n-1)/2$
edges of the complete graph on $[n]$ at least $\lambda$ times,
and since $n_i \leq n$ for all $i$,
it follows that
$$
\frac{n}{4} \sum_{i=1}^m n_i
\geq \sum_{i=1}^m \frac{n_i^2}{4} \geq \lambda n(n-1)/2.
$$
This implies that $\mbox{cap}(\HH)=\sum_{i=1}^m n_i \geq 2\lambda
(n-1)$,
as needed.

Note that for any even $n$ this inequality is tight for
infinitely many (large) values of $\lambda$. In particular it is
tight for $\lambda=\frac{n}{4n-4} {n \choose {n/2}}$, and if there
is a Hadamard matrix of order $n$ then it is tight for $\lambda=n/2$ as
well. In addition, if for some fixed $n$ it is tight for
$\lambda_1$ and $\lambda_2$ then it is also tight for their
sum $\lambda=\lambda_1+\lambda_2$.

We next prove the second inequality, that
$$\mbox{cap}(\HH) \geq 
n [\log n+ \lfloor (\lambda-1)/2 
\rfloor \log (\frac{2 \log n}{\lambda-1})].
$$
Without loss of generality assume that each of the bipartite graphs
$H_i$ in $\HH$ is a complete bipartite graph, and let
$L_i, R_i \subset [n]$ denote its two color classes.
For each
vertex $j \in [n]$ let $A_j$ denote the set of indices 
$i$ for which 
the vertex $j$ belongs to the vertex class $L_i$ of $H_i$ and let
$B_j$ be the set of indices $i$ for which $j \in R_i$.
Let $x_j=|A_j|+|B_j|$ be the total number of bipartite graphs
$H_i$ that contain the vertex $j$. Note that $x_j \geq \lambda$
for each $j$, as any edge incident with $j$ must be covered at
least $\lambda$ times.

Put $r=\lfloor (\lambda-1)/2 \rfloor$ and let $v=(v_1,v_2, \ldots ,
v_m\}$ be a uniform random binary vector of length $m$. 
For each $j$, $1 \leq j \leq n$, let
$E_j$ denote the event that the number of indices $i$ that belong
to $A_j$ for which $v_i =1$  plus the number of indices $i$ that
belong to $B_j$ for which $v_i=0$ is at most $r$. It is clear that
the probability of $E_j$ is exactly
the probability that the binomial random variable
$B(x_j,1/2)$ is at most $r$, which is
$$
p(x_j,1/2)=\frac{\sum_{q=0}^r {x_j \choose q} }{2^{x_j}}.
$$
Note, crucially, that the events $E_j$ are pairwise disjoint. This
is because for every two distinct vertices $j$ and $j'$ there are at
least $\lambda >2r$ indices $i$ for which $j$ and $j'$ belong to
the two distinct vertex classes of $H_i$. Therefore
$$
\sum_{j=1}^n p(x_j,1/2) \leq 1.
$$
The desired lower bound for $\mbox{cap}(\HH) = \sum_{j=1}^n x_j$ 
can be deduced from the last inequality by a convexity argument.
We proceed with the details. Note, first, that for every $x$,
$p(x,r) \geq (x/r)^r 2^{-x}$. Therefore
\begin{equation}
\label{e21}
\sum_{j=1}^n (x_j/r)^r2^{-x_j} \leq 1.
\end{equation}
Recall also that for each $j$, $x_j \geq \lambda \geq 2r+1$.
Consider the function $f(x)=(x/r)^r 2^{-x}$. 
A simple computation shows that for $r=1$ its second derivative is
$(\ln 2) 2^{-x} [x \ln 2-2]$ which is positive for all
$x \geq 2r+1=3$. For $r \geq 2$ the second derivative of $f(x)$ is
$$
(x/r)^{r-2}2^{-x } [ ((x/r) \ln 2 -1)^2 -1/r].
$$
It is not difficult to check that this is positive for all
$x \geq 2r+1$. This shows that $f(x)$ is convex in the relevant
range. Therefore, by (\ref{e21}) together with Jensen's Inequality,
if we denote $x=\mbox{cap}(\HH)=\sum x_j$ we get
$$
n (\frac{x}{nr})^r 2^{-x/n} \leq 1
$$
implying that
$$
x \geq n[\log n+r \log (\frac{x}{nr})].
$$
Since $x/n \geq \log n$ this shows that
$$
x \geq n[\log n+ r \log (\frac{\log n}{r}]
\geq n[\log n + \lfloor (\lambda-1)/2 \rfloor \log 
(\frac{2 \log n}{\lambda-1})].
$$
This completes the proof of Theorem \ref{t11}.  \hfill $\Box$

\subsection{The upper bound}
In this subsection we prove Theorem \ref{t12}.  Put $k=k(n,\lambda)$ and let
$A=(a_{ij})$ be the $k$ by $n$ binary matrix  whose columns are $n$
of the codewords of a binary code of length $k$ with minimum
distance  (at least) $\lambda$. For each $i$, $1 \leq i \leq k$,
let $H_i$ be the complete bipartite graph on the classes of vertices
$L_i=\{j: a_{ij}=0\}$  and $R_i=\{j: a_{ij}=1\}.$
It is easy to see that these bipartite graphs cover every edge of
the complete graph on $[n]$ at least $\lambda$ times. The capacity
of this covering is at most $kn$, establishing the first part
of the theorem. The subsequent items in the theorem follow by
considering appropriate known error correcting codes, see, e.g.
\cite{MS}. 

For the first item simply take the code consisting of all $2^{k-1}$
codewords  with even Hamming weight. Since $2^{k-1} \geq n$ for
$k=\lceil \log n \rceil +1$ the claimed result follows.
The second and third items
follow from the
Gilbert-Varshamov bound which gives that the maximum cardinality of
a binary code with length $k$ and distance $\lambda$ is at least
$$
\frac{2^k}{\sum_{i=0}^{\lambda-1} {k \choose i}}.
$$
This quantity is at least
$2^{k}(\frac{ek}{\lambda-1})^{-(\lambda-1)}$, implying the second
item. For any $\lambda=ck \leq k/2$ this
quantity is also at least $2^{(1-H(c)k}$, where $H(x)$ is the
binary entropy function. This yields
the third item.

The fourth follows by considering an appropriate augmented BCH
code. For any $k$ which is a power of $2$ and for any $d$ this is
a (linear) binary code of length $k$ with 
$$
n=\frac{2^k}{2k^{d-1}}
$$ 
codewords and minimum distance $2d$.  
For $d-1=\lfloor (\lambda-1)/2 \rfloor$, $2d \geq \lambda$ and
$$
k = \log n+1 +(d-1) \log k \leq
\log n + \lfloor (\lambda-1)/2 \rfloor \log \log n+2.
$$
This completes the proof of
Theorem \ref{t12}.  \hfill $\Box$

\section{General graphs}

In this section we prove Theorem \ref{t14}. 
We need the following simple lemma.
\begin{lemma}
\label{l41}
Let $E_i, i \in I$ be a finite collection of events in a (discrete)
probability space. Suppose that for every point $x$ in the space, if
$x \in E_i$ then the total number of events
$E_j$ in the collection that contain $x$ is at most $a_i$. Then
\begin{equation}
\label{e41}
\sum_{i \in I} \frac{\mbox{Prob}(E_i)}{a_i}  \leq 1.
\end{equation}
\end{lemma}
It is worth noting that the above holds (with the same 
proof) for any probability space,
the assumption that it is discrete here is merely because this is the case
we need, and it slightly simplifies the notation in the proof.
\begin{proof}
Let $x$ be an arbitrary point of the space, and let $p(x)$ denote its 
probability. Suppose it belongs to
$r$ of the events $E_i$, let these be
$E_{i_1}, \ldots E_{i_r}$. By the definition of the numbers $a_i$
it follows that $a_{i_j} \geq r$ for all $1 \leq j \leq r$. Therefore
the total contribution of the point $x$ to the sum in the
left-hand-side of (\ref{e41}) is 
$$
\sum_{j=1}^r \frac{p(x)}{a_{i_j}} \leq 
\sum_{j=1}^r \frac{p(x)}{r} \leq  p(x).
$$
The desired result follows by summing over all points $x$ in the space.
\end{proof}
\vspace{0.2cm}

\noindent
{\bf Proof of Theorem \ref{t14}}:\, 
Let $G$ be a graph on the set of vertices $[n]$, let $\alpha_i$
denote the maximum cardinality of an independent set of $G$ containing
the vertex $i$, and
let $\HH=\{H_1, H_2, \ldots ,H_m\}$ be a bipartite covering of $G$.

As in the proof of Theorem \ref{t11} we may and will assume, without
loss of generality,  that each of the bipartite graphs
$H_i$ in $\HH$ is a complete bipartite graph. Let
$L_i, R_i \subset [n]$ denote its two color classes.
For each
vertex $j \in [n]$ let $A_j$ denote the set of indices 
$i$ for which 
the vertex $j$ belongs to the vertex class $L_i$ of $H_i$ and let
$B_j$ be the set of indices $i$ for which $j \in R_i$.
Let $x_j=|A_j|+|B_j|$ be the total number of bipartite graphs
$H_i$ that contain the vertex $j$. Our objective is to prove a lower
bound for the capacity of $\HH$, which is exactly the sum
$\sum_{j=1}^n x_j$.

Let $v=(v_1,v_2, \ldots , v_m\}$ be a uniform random binary 
vector of length $m$. 
For each $j$, $1 \leq j \leq n$, let
$E_j$ denote the event that $v_i=0$ for every index $i$ that belongs
to $A_j$ and $v_i=1$ for every index $i$ that belongs to 
$B_j$. Note that the probability of $E_j$ is exactly
$2^{-x_j}$. Note also that if some point $v=(v_1,v_2, \ldots , v_m\}$
belongs  to the events $E_j, j \in J$, then the set of vertices
$J \subset [n]$ is an independent set of $G$. Indeed, if some
two vertices in $J$ are adjacent, then the edge connecting them
belongs to at least one of the graphs $H_i$ implying that one of these
vertices belongs to $L_i$ whereas the other lies in $R_i$ and showing
that they can't both satisfy the requirement given by $v_i$. 
It thus follows that any point $v$ that lies in $E_j$ belongs to
at most $\alpha_j$ of the events $E_{j'}$. Therefore, by Lemma
\ref{l41}
$$
\sum_{j=1}^n 2^{-x_j-\log \alpha_j}= \sum_{j=1}^n \frac{2^{-x_j}}{\alpha_j}
\leq 1.
$$
By the arithmetic-geometric means inequality this implies
$$
n2^{-(\sum_{j=1}^n x_j +\sum_{j=1}^n \log \alpha_j)/n} \leq 1,
$$
giving 
$$
2^{\sum_{j=1}^n x_j }  
\geq n^n2^{-\sum_{j=1}^n \log \alpha_j}=
2^{\sum_{j=1}^n (\log n-\log \alpha_j)}.
$$
Therefore
$$
\sum_{j=1}^n x_j \geq \sum_{j=1}^n \log (\frac{n}{\alpha_j}),
$$
completing the proof.   \hfill $\Box$

\end{document}